\documentclass{article}

\usepackage{amssymb}
\usepackage{amsmath}
\usepackage{mathrsfs}
\usepackage{cite}
\usepackage{amsthm}
\usepackage{enumerate}
\theoremstyle{definition}
\numberwithin{equation}{section}

\newtheorem{lemma}{Lemma}[section]
\newtheorem{theorem}{Theorem}[section]

\newtheorem{corollary}{Corollary}[section]

\begin{document}

\setlength\abovedisplayskip{2pt}
\setlength\abovedisplayshortskip{0pt}
\setlength\belowdisplayskip{2pt}
\setlength\belowdisplayshortskip{0pt}

\title{\bf Dzhaparidze-van Zanten type inequalities for self-normalized martingales \author{WenCong Zhang}} \maketitle
 \footnote{Received: \today.}
 \footnote{MR Subject Classification: 0211 .}
\begin{center}
\begin{minipage}{135mm}
{\bf \small Abstract}.\hskip 2mm {\small
The Bernstein inequality is a tight upper bound on tail probabilities for independent random variables. Freedman extended the Bernstein inequality to martingales with differences bounded from above, and then Dzhaparidze and van Zanten generalized Freedman's result to non-bounded locally square integrable martingales. In this paper, we derive some Dzhaparidze-van Zanten type inequalities for self-normalized martingales with square integrable and non-square integrable differences.}
\end{minipage}
\end{center}
\section{Introduction}
The classical Bernstein inequality gives a tight upper bound on tail probabilities for sums of independent   random variables.
 Let $\{\xi_i\}_{i=0,1,\cdots}$ be a
sequence of zero-mean independent random variables satisfying  $|\xi_i|\leq a$ for all $i$ and some positive constant $a$.
 Denote  $S_n= \sum_{i=1}^n\xi_i$ the partial sum of $\xi_i.$ The Bernstein inequality implies the following result: for all $z>0$,
\begin{equation}
\mathbb{P}( S_n >z)\leq  \exp\big\{-\frac{1}{2}\big(\frac{z^2}{var(S_n)+\frac{az}{3}}\big)\big\},
\end{equation}
where  $var(S_n)=\sum_{i=1}^n\mathbb{E}\xi_i^2$ is the variance of $S_n$.

Freedman \cite{F75} extended Bernstein's result to the case of discrete-time martingales. Let $\{\xi_i, \mathfrak{F}_n\}_{n=1,2,\cdots}$ be a sequence of  martingale differences satisfying $|\xi_i|\leq a$ for some constant $a$.
Then, by definition,  $\{S_n, \mathfrak{F}_n\}_{n=1,2,\cdots}$   is a martingale. Denote by $\langle S\rangle_n$   the conditional variance of $S_n$, that is
\[\langle S\rangle_n=\sum_{i=1}^{n}\mathbb{E}[\xi_i^2|\mathfrak{F}_{i-1}].\]
The Freedman  inequality states that for all $x,L\geq 0$,
\begin{equation}\label{eq1.2}
\mathbb{P}\big( S_n >x, \langle S\rangle_n \leq L\big)\leq  \exp\big\{-\frac{1}{2}\big(\frac{x^2}{L+\frac{ax}{3}}\big)\big\}.
\end{equation}
The Freedman  inequality is further generalized by Dzhaparidze and van Zanten \cite{DV01} to
martingales with non-bounded differences.   Define the second-order process $\{H_n^a\}_{n=1, 2,...}$ as follows
\[H_n^a =\sum_{i=1}^{n}\xi_i^21_{\{|\xi_i|>a\}}+\langle S\rangle_n.\]
Then for every $x, L>0$,
\begin{equation}\label{eq1.3}
\mathbb{P}\big( S_n >x, H_n^a\leq L\big)\leq  \exp\big\{-\frac{1}{2}\big(\frac{x^2}{L}\big)\psi(\frac{ax}{L})\big\},
\end{equation}
where $\psi$ is defined by
\[\psi(x)=\frac{2}{x^2}\int_{0}^{x}\log(1+y)dy.\]
and satisfies
\[\psi(x)\geq \frac{1}{1+x/3}\quad \textrm{if}\ \ x\geq -1.\]
It is obviously that the new bound in inequality (\ref{eq1.3}) is somewhat sharper than the earlier bound in (\ref{eq1.2}). Besides, under the earlier condition that $|\xi_i|\leq a$ for all $i$, the first term in $H_n^a$ vanishes, then $H_n^a=\langle S\rangle_n$. So inequality (\ref{eq1.3})   implies the Freedman inequality and, as well as its consequence, the classical Bernstein inequality.

Despite the Bernstein inequality for martingale is well studied, there are only a few of results on Bernstein type inequalities for self-normalized martingales.
When the martingale differences are conditional symmetric in distribution with respect to zero, de la Pe\~{n}a \cite{D99} have established the following inequalities for self-normalized martingales: for all $x>0$,
\begin{equation}\label{eq1.4}
\mathbb{P}\big(\frac{S_n}{[S]_n}\geq x\big)\leq \sqrt{\mathbb{E}\big[\exp\big\{-\frac{1}{2}x^2[S]_n\big\}\big]},
\end{equation}
and for all $x,y>0$,
\begin{equation}\label{eq1.5}
\mathbb{P}\big(\frac{S_n}{[S]_n}\geq x,[S]_n\geq y\big)\leq \exp\big\{-\frac{1}{2}x^2y\big\},
\end{equation}
where $[S]_n=\sum_{i=1}^{n}\xi_i^2$.
de la Pena and Pang \cite{DP09} generalized inequality (\ref{eq1.4}) to self-normalized processes. Let $(A,B)$ be a pair of random variables with $B>0$ in the probability space $(\Omega,\mathfrak{F},P)$ satisfies the canonical assumption
\[\mathbb{E}\big[\exp\big\{\lambda A-\frac{\lambda^2}{2}B^2\big\}\big]\leq 1,\qquad\lambda \in \mathbb{R}.\]
Suppose $\mathbb{E}[|A|^p]<\infty$ for some $p\geq 1$. Then for any $x>0$ and for $q\geq 1$ such that $1/p+1/q=1$,
\begin{equation}\label{eq1.6}
\mathbb{P}\Big(\frac{|A|}{\sqrt{\frac{2q-1}{q}(B^2+(\mathbb{E}[|A|^p])^{2/p})}}\geq x\Big)\leq \big(\frac{q}{2q-1}\big)^{\frac{q}{2q-1}}x^{-\frac{q}{2q-1}}\exp\big\{-\frac{x^2}{2}\big\}.
\end{equation}

Recall that an integrable random variables $X$ is called heavy on left if $\mathbb{E}X=0$ and, for all $a>0$, $\mathbb{E}[T_a(x)]\leq 0$, where
\[T_a(x)=\min(|X|,a)\textrm{sign}(X)\]
is the truncated version of $X$.  Bercu and Tuati \cite{BT08} extended (\ref{eq1.4}) and (\ref{eq1.5}) when $S_n$ is a locally square integrable martingale heavy on left: for all $x>0,a\geq 0$ and $b>0$,
\begin{equation}\label{eq1.7}
\mathbb{P}\big(\frac{S_n}{a+b[S]_n}\geq x\big)\leq \inf_{p>1}\big(\mathbb{E}\big[\exp\big\{-(p-1)x^2\big(ab+\frac{b^2}{2}[S]_n\big)\big\}\big]\big)^{1/p}
\end{equation}
and, for all $y>0$,
\begin{equation}\label{eq1.8}
\mathbb{P}\big(\frac{S_n}{a+b[S]_n}\geq x,[S]_n\geq y\big)\leq \exp\big(-x^2\big(ab+\frac{b^2y}{2}\big)\big).
\end{equation}
Obviously that letting $a=0$ and $b=1$ in (\ref{eq1.7}) and (\ref{eq1.8}), we can respectively get (\ref{eq1.4}) and (\ref{eq1.5}).

Recently,  Fan and Wang \cite{FW19} obtained some results similar to inequalities (\ref{eq1.4}) and (\ref{eq1.5}) for self-normalized martingales with differences bounded from below.

In this paper, we aim to extend inequality (\ref{eq1.3}) to  self-normalized martingales,  with   Dzhaparidze-van Zanten type self-normalized factors. Based on inequality (\ref{eq1.8}), we obtain another inequality for self-normalized martingales with the condition of heavy on left.

The paper is organized as follows. We present our main results in Section 2. In Section 3, we discuss some applications of our results, and the proofs of our main results are given in Section 4.
\section{Main Result}
Let $(\xi_i,\mathfrak{F}_i)_{i=0,1,\cdots }$ be a  sequence of real-valued  martingale differences defined on a probability space $(\Omega,\mathfrak{F},P)$, where $\xi_0=0$ and $\{\emptyset,\Omega\}=\mathfrak{F}_0\subseteq\cdots\subseteq \mathfrak{F}_n\subseteq\mathfrak{F}$ are increasing $\sigma$-fields. So by definition, we have $\mathbb{E}[\xi_i|\mathfrak{F}_{i-1}]=0,i=1,2,\cdots$. Set
\[S_0=0\qquad \textrm{and}\qquad S_k=\sum_{i=1}^{k}\xi_i\]
for $k=1,2,\cdots,n$. Then $S=(S_k,\mathfrak{F}_k)_{k=1,2,\cdots }$ is a martingale.

First,  we consider the martingales with squared integrable  differences.
Given $y\geq0.$ Let $[S]_n(y)$ and $\langle S\rangle_n(y)$ respectively be
\[[S]_n(y)=\sum_{i=1}^{n}(\xi_i)^21_{\{\xi_i>y\}}\qquad \textrm{and}\qquad \langle S\rangle_n(y)=\sum_{i=1}^{n}\mathbb{E}[\xi_i^21_{\{\xi_i\leq y\}}|\mathfrak{F}_{i-1}].\]
Denote
\[
B_n(y)=[S]_n(y)+\langle S\rangle_n(y).
\]
Our first result is the following Dzhaparidze-van Zanten type inequalities for self-normalized martingales.
\begin{theorem}\label{th2.1}
For all $x, y\geq 0$,
\begin{eqnarray}
 \mathbb{P}(\frac{S_n}{B_n(y)} \geq x)
&\leq&\inf_{p>1}\big(\mathbb{E}\big[\exp\big\{-(p-1)f(x,y)B_n(y)\big\}1_{\{S_n\geq xB_n(y)\}}\big]\big)^{1/p} \label{ineq2.1}\\
&\leq& \inf_{p>1}\big(\mathbb{E}\big[\exp\big\{-(p-1)\frac{x^2}{2(1+{xy}/3)}B_n(y)\big\}1_{\{S_n\geq xB_n(y)\}}\big]\big)^{1/p},  \label{ineq2.2}
\end{eqnarray}
where
\[
f(x,y)=\frac{xy(\ln (xy+1)-1)+\ln (xy+1)}{y^2}
\]
with the convention that $f(x,0)=\frac{x^2}{2}$. Inequality (\ref{ineq2.2}) implies that for all $x, y\geq 0$,
\begin{equation}\label{ineq2.3}
\mathbb{P}(\frac{S_n}{B_n(y)}\geq x)  \leq\inf_{p>1}\big(\mathbb{E}\big[\exp\big\{-(p-1)\frac{x^2}{2(1+{xy}/3)}B_n(y)\big\}\big]\big)^{1/p}.
\end{equation}
Moreover, we also have for all$x,y,z>0$,
\begin{equation}
\mathbb{P}(\frac{S_n}{B_n(y)}\geq x,B_n(y)\leq z) \leq \exp\big\{-\frac{x^2z}{2(1+\frac{xy}{3})}\big\}.
\end{equation}
\end{theorem}

Notice that $H_n^y\leq B_n(y)$. Inspiring the proof of Theorem \ref{th2.1}, it is easy to see that the inequalities (\ref{ineq2.1})-(\ref{ineq2.3}) hold also when
$B_n(y)$ is replaced by  $H_n^y$.

 Let $[S]_n^{+}$ and $\langle S\rangle_n^{-}$ be, respectively, the positive term of the squared variance and the negative term of the conditional variance of the martingale $S_n$, that is
\[[S]_n^{+}=\sum_{i=1}^{n}(\xi_i^{+})^2\qquad \textrm{and}\qquad \langle S\rangle_n^{-}=\sum_{i=1}^{n}\mathbb{E}[(\xi_i^{-})^2|\mathfrak{F}_{i-1}],\]
where
\[\xi_i^{+}=\max\{\xi_i,0\}\qquad\textrm{and}\qquad\xi_i^{-}=\max\{-\xi_i,0\}.\]
Clearly, it holds $ B_n(0) =[S]_n^{+}+\langle S\rangle_n^{-}$.
Taking $y=0$ in Theorem \ref{th2.1}, we have the following corollary.
\begin{corollary}\label{c2.1}
For all $x>0$,
\begin{equation}
\begin{aligned}
&\mathbb{P}(\frac{S_n}{B_n(0)}\geq x)\\
&\leq \inf_{p>1}\big(\mathbb{E}[\exp\{-(p-1)\frac{x^2}{2}(B_n(0))\}1_{\{S_n\geq x(B_n(0))\}}]\big)^{1/p} \nonumber \\
&\leq \inf_{p>1}\big(\mathbb{E}[\exp\{-(p-1)\frac{x^2}{2}(B_n(0))\}]\big)^{1/p}
\end{aligned}
\end{equation}
and, for all $x, y>0$,
\begin{equation}
\mathbb{P}(\frac{S_n}{B_n(0)}\geq x,\, B_n(0) \geq y)\leq \exp\{-\frac{x^2y}{2}\}.
\end{equation}

The last inequality can be regarded as a self-normalized version of Delyon's inequality \cite{De2009} , where Delyon proved that for all $x, y>0$,
\begin{equation}
\mathbb{P}( S_n \geq x,\, B_n(0)\leq y)\leq \exp\{-\frac{x^2}{2y}\}.
\end{equation}

When the normalized factor $B_n(y)$ in Theorem \ref{th2.1} is replaced by $\sqrt{B_n(y)}$. We have the following results.
\begin{theorem}\label{th2.2}
For all $b>0$, $M\geq 1$  and $x, y\geq 0$,
\begin{equation}
\begin{aligned}
\mathbb{P}(\frac{S_n}{\sqrt{B_n(y)}}\geq x,&b\leq \sqrt{B_n(y)}\leq bM) \\
&\leq \sqrt{e}(1+2(1+x)\ln M)\exp\big\{-\frac{x^2}{2(1+\frac{xy}{3b})}\big\}.
\end{aligned}
\end{equation}
\end{theorem}

Taking $y=0$ in Theorem \ref{th2.2}, we have the following corollary.
\end{corollary}
\begin{corollary}\label{col2.2}
For all $b>0$, $M\geq 1$ and $x>0$,
\begin{equation}\label{eq2.7}
\begin{aligned}
\mathbb{P}(\frac{S_n}{\sqrt{B_n(0)}}\geq x,\ &b\leq \sqrt{B_n(0)}\leq bM)\\
&\leq \sqrt{e}(1+2(1+x)\ln M)\exp\{-\frac{1}{2}x^2\}.
\end{aligned}
\end{equation}
\end{corollary}

Next we consider the case that $\xi_i$ are heavy on left.

\begin{theorem}\label{th2.5}
Assume that $\xi_i$ are heavy on left for all $i$. Then for all $b>0$, $M\geq 1$ and $x>0$,
\begin{equation}
\begin{aligned}
\mathbb{P}(\frac{S_n}{\sqrt{[S]_n}}\geq x,&b\leq \sqrt{[S]_n}\leq bM)\\
&\leq \sqrt{e}(1+2(1+x)\ln M)\exp\{-\frac{1}{2}x^2\}.
\end{aligned}
\end{equation}
\end{theorem}
Now we consider the martingales with non-squared-integrable differences. Denote
\[[S]_n^{+}(\beta)=\sum_{i=1}^{n}(\xi_i^{+})^\beta\qquad \textrm{and}\qquad \langle S\rangle_n^{-}(\beta)=\sum_{i=1}^{n}\mathbb{E}[(\xi_i^{-})^\beta|\mathfrak{F}_{i-1}].\]
Denote
\[G_n(\beta)=[S]_n^{+}(\beta)+\langle S\rangle_n^{-}(\beta).\]
We have the following inequalities for self-normalized martingales.
\begin{theorem}\label{th2.3}
If $\mathbb{E}|\xi_i|^\beta<\infty$ for some $\beta \in (1,2)$, then for all   $x>0$,
\begin{equation}
\begin{aligned}
&\mathbb{P}(\frac{S_n}{G_n(\beta)}\geq x)\\
&\leq\inf_{p>1}\big(\mathbb{E}\big[\exp\big\{(p-1)(\frac{x}{\beta})^{\frac{\beta}{\beta-1}}(1-\beta)G_n(\beta)\big\}1_{\{S_n\geq xG_n(\beta)\}}\big]\big)^{1/p}\\
&\leq\inf_{p>1}\big(\mathbb{E}\big[\exp\big\{(p-1)(\frac{x}{\beta})^{\frac{\beta}{\beta-1}}(1-\beta)G_n(\beta)\big\}\big]\big)^{1/p}.
\end{aligned}
\end{equation}
\end{theorem}

\begin{theorem}\label{th2.4}
If $\mathbb{E}|\xi_i|^\beta<\infty$ for some $\beta \in (1,2)$, then for all $b>0$, $M\geq 1$ and $x>0$,
\begin{equation}
\begin{aligned}
\mathbb{P}\big(\frac{S_n}{\sqrt[\beta]{G_n(\beta)}}&\geq x, b^{\frac{1}{\beta-1}}\leq \sqrt[\beta]{G_n(\beta)}\leq (bM)^{\frac{1}{\beta-1}}\big)\\
&\leq (1+2(1+x)\ln M)\exp\big\{-(\frac{x}{\beta})^{\frac{\beta}{\beta-1}}(1-\frac{1}{\beta})\big\}.
\end{aligned}
\end{equation}
\end{theorem}
\section{Application}
\subsection{Student's $t$-statistic}
Recall that Student's $t$-statistic $T_n$ is defined by
\[T_n=\frac{\sqrt{n}\bar{\xi}}{\big(\frac{1}{n-1}\sum_{j=1}^{n}(\xi_j-\bar{\xi})^2\big)^{1/2}},\]
where $\bar{\xi}=\frac1n\sum_{i=1}^{n}\xi_i$. Clearly, the following equation is true:
\[T_n=\frac{S_n}{\sqrt{[S]_n}}\big(\frac{n-1}{n-(S_n/\sqrt{[S]_n})^2}\big)^{1/2}.\]
Since $x/(n-x^2)^{1/2}$ is increasing in $(-\sqrt{n},\sqrt{n})$, so
\begin{equation}\label{eq2.20}
\{T_n\geq x\}=\big\{\frac{S_n}{\sqrt{[S]_n}}\geq x\big(\frac{n}{n+x^2-1}\big)^{1/2}\}.
\end{equation}
 Using Theorem \ref{th2.5} and equation (\ref{eq2.20}), we get the following deviation inequality for $t$-statistic.
\begin{theorem}
Assume that $\xi_i$ are heavy on left for all $i$. Then for all $b>0$, $M\geq 1$ and $x>0$,
\begin{equation}
\begin{aligned}
&\mathbb{P}(T_n\geq x, b\leq \sqrt{[S]_n}\leq bM)\\
&\leq \sqrt{e}\big(1+2\big(1+x\big(\frac{n}{n+x^2-1}\big)^{1/2}\big)\ln M\big)\exp\{-\frac{nx^2}{2(n+x^2-1)}\}.
\end{aligned}
\end{equation}
\end{theorem}
\subsection{Linear Regressions}

The stochastic linear regression can be expressed for all $n\geq 0$:
\[X_{n+1}=\theta \phi_n+\varepsilon_{n+1},\]
where $X_n,\phi_n,\varepsilon_n$ here are respectively called the observation,the regression variable and the drive noise. We assume that $\phi_n$ is a sequence of i.i.d. random variables and $\varepsilon_n$ is a sequence of identically distributed random variables with mean zero and variance $\sigma^2>0$. Moreover, we suppose that for all $n\geq 0$, the random variable $\varepsilon_{n+1}$ is independent of $\mathfrak{F}_n$ where $\mathfrak{F}_n=\sigma(\phi_0,\epsilon_1,\cdots,\phi_{n-1},\varepsilon_n)$. We give the least-squares estimator $\hat{\theta}_n$ as, for all $n\geq 1$,
\[\hat{\theta}_n=\frac{\sum_{k=1}^{n}\phi_{k-1}X_k}{\sum_{k=1}^{n}\phi_{k-1}^2}.\]
Bercu and Touati \cite{Be2018} has given the convergence rate of $\hat{\theta}_n-\theta$ when $(\phi_n)$ and $(\varepsilon_n)$ are normal random variables. Now, we would like to give a a convergence rate of  $\hat{\theta}_n-\theta$ when $(\varepsilon_n)$ only has the upper bound using Theorem \ref{th2.1}.
\begin{theorem}\label{th3.2}
Assume that $\varepsilon_i\leq y$ for all $y>0$ and all $i$. If $|\phi_i|\leq 1$, then for all $x>0$,
\begin{equation}
\begin{aligned}
\mathbb{P}=(|\hat{\theta}_n-\theta|\geq x)&\leq \\
&2\inf_{p>1}\big(\mathbb{E}\big[\exp\big\{-(p-1)\frac{x^2}{2(\sigma^2+\frac{xy}{3}))}\sum_{k=1}^{n}\phi_{k-1}^2\big\}\big]\big)^{1/p}.
\end{aligned}
\end{equation}
\end{theorem}
By Theorem \ref{th2.2}, we obtain the following result.
\begin{theorem}\label{th3.3}
Assume that $\varepsilon_i\leq y$ for all $y>0$ and all $i$. If $|\phi_i|\leq 1$, then for all $b>0, M\geq 1$ and $x>0$,
\begin{equation}
\begin{aligned}
\mathbb{P}\big(|\hat{\theta}_n-\theta|\sqrt{\sum_{k=1}^{n}\phi_{k-1}^2}\geq x,b\leq \sqrt{\sum_{k=1}^{n}\phi_{k-1}^2}\leq bM\big)&\leq \\
&2\sqrt{e}\big(1+2(1+\frac{x}{\sigma})\ln M\big)\exp\big\{-\frac{x^2}{2(\sigma^2+ xy/3b)}\big\}.
\end{aligned}
\end{equation}
\end{theorem}
\subsection{Stochastic TSP Problem}
In the stochastic modeling of the TSP, let $X_1,X_2,\cdots,X_n$ be i.i.d. unifromly distributed on $[0,1]^d(d\geq 2)$ and $T_n$ be the shortest closed path through the $n$ random points $X_1,X_2,\cdots,X_n$. In particular, by Azuma's inequality (see Theorem 2.1 in Steele \cite{ste81}), we have for $n\geq 2$ and some constant $C$,
\begin{equation}
\mathbb{P}(|T_n-\mathbb{E}[T_n]|\geq t)\leq \left\{
\begin{aligned}
&\exp\{-t^2/(C\log n)\},& \text{for}\quad d=2,\\
&\exp\{-t^2/(Cn^{(d-2)/d})\},& \text{for}\quad d>2.
\end{aligned}\right.
\end{equation}
In particular, $var(T_n)<\infty$ for all $n$. Cerf et al. \cite{Cerf97} did numerical simulation for $d=2$ to support that the variable $T_n-\mathbb{E}[T_n]$ should have a Gaussian distribution as $n\to \infty$. Here we obtain the following upper bound.
\begin{theorem}\label{th3.4}
For the stochastic TSP problem,
\begin{equation}
\begin{aligned}
\mathbb{P}\big(\frac{T_n-\mathbb{E}(T_n)}{\sqrt{\sum_{i=1}^{n}d_i^2}}\geq t, &c_1(d)n^{1/2-1/d}\leq \sqrt{\sum_{i=1}^{n}d_i^2}\leq c_1(d)n^{1/2}\big)\leq\\
&\sqrt{e}(1+\frac{2}{d}(1+t)\ln n)\exp\{-\frac{t^2}{2}\},
\end{aligned}
\end{equation}
where $t>0$ and $c_1(d)>0$ and $d_i=\mathbb{E}[T_n|\mathfrak{F}_i]-\mathbb{E}[T_n|\mathfrak{F}_{i-1}]$ for $1\leq i\leq n$, and $\mathfrak{F}_i=\sigma\{X_1,\cdots,X_i\}$.
\end{theorem}
\section{Proof of Theorems}

\subsection{Proof of Theorem \ref{th2.1}}
To prove Theorem \ref{th2.1}, we need the following lemma (see  Lemma 5.1  of Fan et al.\cite{FGL17}).
\begin{lemma}\label{t1.3}
For all $\lambda, y\geq 0$, denote
\[U_n(\lambda)=\exp\big\{\lambda S_n-\big(\frac{e^{\lambda y}-1-\lambda y}{y^2}\big)B_n(y)\big\}.\]
Then $(U_i(\lambda),\mathfrak{F}_i)_{i=0,\cdots,n}$ is a supermartingale and satisfies that
\[\mathbb{E}[U_n(\lambda)]\leq 1.\]
\end{lemma}
We use the method of Bercu and Touati \cite{BT08}. Let $W_n=\{S_n\geq xB_n(y)\},x>0,y\geq 0$. By Markov's inequality, H\"{o}lder inequality and Lemma \ref{t1.3}, we have for all $\lambda>0$ and $q>1$,
\begin{equation*}
\begin{aligned}
\mathbb{P}(W_n)&\leq \mathbb{E}\big[\exp\big\{\frac{\lambda}{q}\big(S_n-xB_n(y)\big)\big\}1_{W_n}\big]\\
&=\mathbb{E}\big[\exp\big\{\frac{1}{q}\big(\lambda S_n-\frac{e^{\lambda y}-1-\lambda y}{y^2}B_n(y)\big)\big\}\\
&\exp\big\{\frac{1}{q}\big(\frac{e^{\lambda y}-1-\lambda y}{y^2}-\lambda x\big)B_n(y)\big\}1_{W_n}\big]\\
&\leq \big(\mathbb{E}\big[\exp\big\{\frac{p}{q}\big(\frac{e^{\lambda y}-1-\lambda y}{y^2}-\lambda x\big)B_n(y)\big\}1_{W_n}\big]\big)^{1/p}(\mathbb{E}[U_n(\lambda)])^{1/q}\\
&\leq \big(\mathbb{E}\big[\exp\big\{\frac{p}{q}\big(\frac{e^{\lambda y}-1-\lambda y}{y^2}-\lambda x\big)B_n(y)\big\}1_{W_n}\big]\big)^{1/p},
\end{aligned}
\end{equation*}
where $p=1+p/q$. Consequently, as $p/q=p-1$, we   deduce that
\[\mathbb{P}(W_n)\leq \inf_{p>1}\big(\mathbb{E}\big[\exp\big\{-(p-1)\big(\lambda x-\frac{e^{\lambda y}-1-\lambda y}{y^2}\big)B_n(y)\big\}1_{W_n}\big]\big)^{1/p}.\]
The right-hand side of the last inequality attains its minimum at
\[\lambda^{\ast}=\frac{\ln(xy+1)}{y},\]
therefore we obtain
\[\mathbb{P}(W_n)\leq \inf_{p>1}\big(\mathbb{E}\big[\exp\big\{-(p-1)f(x,y)B_n(y)\big\}1_{W_n}\big]\big)^{1/p}.\]
Using the inequality
\[e^{\lambda y}-1-\lambda y\leq \frac{\lambda^2y^2}{2(1-{\lambda y}/3)},\quad\lambda \in [0,\frac{3}{y}),\]
we get for all $x, y\geq0$,
\begin{equation*}
\begin{aligned}
f(x,y)&=\inf_{\lambda\geq 0}\big(\frac{e^{\lambda y}-1-\lambda y}{y^2}-\lambda x\big)\\
&\leq \inf_{\lambda\geq 0}\big(\frac{\lambda^2}{2(1-{\lambda y}/3)}-\lambda x\big)\\
&\leq -\frac{x^2}{2(1+{xy}/3)}.
\end{aligned}
\end{equation*}
Thus, we obtain for all $x, y\geq0$,
\begin{equation*}
\begin{aligned}
\mathbb{P}(W_n)&\leq \inf_{p>1}\big(\mathbb{E}\big[\exp\big\{-(p-1)\big(\lambda x-\frac{e^{\lambda y}-1-\lambda y}{y^2}\big)B_n(y)\big\}1_{W_n}\big]\big)^{1/p}\\
&\leq \inf_{p>1}\big(\mathbb{E}\big[\exp\big\{-(p-1)\frac{x^2}{2(1+{xy}/3)}B_n(y)\big\}1_{W_n}\big]\big)^{1/p},
\end{aligned}
\end{equation*}
which gives the desired inequalities.
\subsection{Proof of Theorem \ref{th2.2}}
The proof of Theorem \ref{th2.2} is based on a modified method of Lipster and Spokoiny \cite{LS00}. Given $a>1$, introduce the geometric series $b_k=ba^k$ and define random events
\[Q_k=\big\{\frac{S_n}{\sqrt{B_n(y)}}\geq x,b_k\leq \sqrt{B_n(y)}\leq b_{k+1}\big\},\quad k=0,1,\cdots,K,\]
where $K$ stands for the integer part of $\log_aM$. Clearly, it holds
\[\big\{\frac{S_n}{\sqrt{B_n(y)}}\geq x,b\leq \sqrt{B_n(y)}\leq bM\big\}\subseteq\cup_{k=0}^{K}Q_k,\]
which leads to
\[\mathbb{P}\big(\frac{S_n}{\sqrt{B_n(y)}}\geq x,b\leq \sqrt{B_n(y)}\leq bM\big)\leq \sum_{k=0}^{K}\mathbb{P}(Q_k).\]
Notice that
\[e^{\lambda y}-1-\lambda y\leq \frac{\lambda^2y^2}{2(1-{\lambda y}/3)},\lambda \in [0,\frac{3}{y}).\]
For any $\lambda \in [0,3)$, the last inequality and Lemma \ref{t1.3} together implies that
\[\mathbb{E}\big[\exp\big\{\lambda S_n-\frac{\lambda^2}{2(1-{\lambda y}/3)}B_n(y)\big\}1_{Q_k}\big]\leq 1.\]
Next, taking $\lambda_k=\frac{x}{b_k+{xy}/3}$, for any $x>0$ and $y\geq0$, we obtain
\begin{equation*}
\begin{aligned}
1&\geq \mathbb{E}\big[\exp\big\{\frac{x}{b_k+{xy}/3}S_n-\frac{x^2}{2b_k(b_k+{xy}/3)}B_n(y)\big\}1_{Q_k}\big]\\
&\geq \mathbb{E}\big[\exp\big\{\frac{x^2}{b_k+{xy}/3}\sqrt{B_n(y)}-\frac{x^2}{2b_k(b_k+{xy}/3)}B_n(y)\big\}1_{Q_k}\big]\\
&\geq \mathbb{E}\big[\exp\big\{\inf_{b_k<c<b_{k+1}}\big(\frac{x^2c}{b_k+{xy}/3}-\frac{x^2c^2}{2b_k(b_k+{xy}/3)}\big)\big\}1_{Q_k}\big]\\
&=\mathbb{E}\big[\exp\big\{\frac{x^2b_{k+1}}{b_k+{xy}/3}-\frac{x^2b_{k+1}^2}{2b_k(b_k+{xy}/3)}\big\}1_{Q_k}\big],
\end{aligned}
\end{equation*}
which implies that
\begin{equation*}
\begin{aligned}
\mathbb{P}(Q_k)&\leq \exp\big\{\frac{x^2b_{k+1}^2}{2b_k(b_k+{xy}/3)}-\frac{x^2b_{k+1}}{b_k+{xy}/3}\big\}\\
&\leq \exp\big\{-\frac{x^2}{1+\frac{xy}{3b_k}}(a-\frac{a^2}{2})\big\}\\
&\leq \exp\big\{-\frac{x^2}{1+\frac{xy}{3b}}(a-\frac{a^2}{2})\big\}.
\end{aligned}
\end{equation*}
Finally, we may pick $a$ to make the right-hand side of the last bound as small as possible. This leads to the choice $a=1+1/(1+x)$, so that
\[x^2(a-\frac{a^2}{2})=x^2\big(1+\frac{1}{1+x}-\frac{1}{2}(1+\frac{1}{1+x})^2\big)\geq \frac{1}{2}(x^2-1).\]
Since $\log(1+1/(1+x))\geq 1/(2(1+x))$ for $x\geq 0$, we obtain $\log_aM\leq 2(1+x)\ln M$, which gives the desired inequality.

\subsection{Proof of Theorem \ref{th2.5}}
The proof of Theorem \ref{th2.5} is similar with the proof of Theorem \ref{th2.2}. Given $a>1$, introduce the geometric series $b_k=ba^k$ and define random events
\[C_k=\big\{\frac{S_n}{\sqrt{[S]_n}}\geq x,b_k\leq \sqrt{[S]_n}\leq b_{k+1}\big\},\quad k=0,1,\cdots,K,\]
where $K$ stands for the integer part of $\log_aM$. Clearly, it holds
\[\big\{\frac{S_n}{\sqrt{[S]_n}}\geq x,b\leq \sqrt{[S]_n}\leq bM\big\}\subseteq \cup_{k=0}^{K}C_k,\]
which leads to
\[\mathbb{P}\big(\frac{S_n}{\sqrt{[S]_n}}\geq x,b\leq \sqrt{[S]_n}\leq bM\big)\leq \sum_{k=0}^{K}\mathbb{P}(C_k).\]
From Lemma 3.1 in Bercu and Touati \cite{BT08}, we can get
\[\mathbb{E}\big[\exp\big\{\lambda S_n-\frac{\lambda^2}{2}[S]_n\big\}1_{C_k}\big]\leq 1,\]
where $\lambda >0$. Next, taking $\lambda_k=\frac{x}{b_k}$, for any $x>0$, we obtain
\begin{equation*}
\begin{aligned}
1&\geq \mathbb{E}\big[\exp\big\{\frac{x}{b_k}S_n-\frac{x^2}{2b_k^2}[S]_n\big\}1_{C_k}\big]\\
&\geq \mathbb{E}\big[\exp\big\{\frac{x^2}{b_k}\sqrt{[S]_n}-\frac{x^2}{2b_k^2}[S]_n\big\}1_{C_k}\big]\\
&\geq \mathbb{E}\big[\exp\big\{\inf_{b_k<c<b_{k+1}}\big(\frac{x^2c}{b_k}-\frac{x^2c^2}{2b_k^2}\big)\big\}1_{C_k}\big]\\
&=\mathbb{E}\big[\exp\big\{\frac{x^2b_{k+1}}{b_k}-\frac{x^2b_{k+1}^2}{2b_k^2}\big\}1_{C_k}\big]\\
&=\mathbb{E}\big[\exp\big\{x^2a-\frac{x^2a^2}{2}\big\}1_{C_k}\big],
\end{aligned}
\end{equation*}
which implies that
\[\mathbb{P}(C_k)\leq \exp\{-x^2(a-\frac{a^2}{2})\}.\]
Let $a=1+1/(1+x)$, so that
\[\mathbb{P}(C_k)\leq \exp\big\{-\frac{1}{2}(x^2-1)\big\}.\]
Then
\[\mathbb{P}(\frac{S_n}{\sqrt{[S]_n}}\geq x,b\leq \sqrt{[S]_n}\leq bM)\leq \sqrt{e}(1+2(1+x)\ln M)\exp\{-\frac{1}{2}x^2\},\]
which gives the desired inequality.
\subsection{Proof of Theorem \ref{th2.3}}
In the proof of Theorem \ref{th2.3}, we make use of the following lemma due to   Fan et al.\cite{FGL17}.
\begin{lemma}\label{t1.4}
Assume $\mathbb{E}|\xi_i|^\beta<\infty$ for some $\beta \in (1,2)$.  Denote
\[V_n(\lambda)=\exp\{\lambda S_n-\lambda^\beta G_n(\beta)\}, \ \ \lambda>0. \]
Then $(V_i(\lambda),\mathfrak{F}_i)_{i=0,\cdots,n}$ is a supermartingale and satisfies
\[\mathbb{E}[V_n(\lambda)]\leq 1.\]
\end{lemma}
The method in this proof is similar with Theorem \ref{th2.1}. Let $D_n=\{S_n\geq xG_n(\beta)\},x>0$. By Markov's inequality, H\"{o}lder inequality and Lemma \ref{t1.4}, we have for all $\lambda>0$ and $q>1$,
\begin{equation*}
\begin{aligned}
\mathbb{P}(D_n)&\leq \mathbb{E}\big[\exp\big\{\frac{\lambda}{q}\big(S_n-xG_n(\beta)\big)\big\}1_{D_n}\big]\\
&=\mathbb{E}\big[\exp\big\{\frac{1}{q}\big(\lambda S_n-\lambda^\beta G_n(\beta)\big)\big\}\\
&\exp\big\{\frac{1}{q}\big(\lambda^\beta-\lambda x\big)G_n(\beta)\big\}1_{D_n}\big]\\
&\leq \big(\mathbb{E}\big[\exp\big\{\frac{p}{q}\big(\lambda^\beta-\lambda x\big)G_n(\beta)\big\}1_{D_n}\big]\big)^{1/p}(\mathbb{E}[V_n(\lambda)])^{1/q}\\
&\leq \big(\mathbb{E}\big[\exp\big\{\frac{p}{q}\big(\lambda^\beta-\lambda x\big)G_n(\beta)\big\}1_{D_n}\big]\big)^{1/p},
\end{aligned}
\end{equation*}
where $p=1+p/q$. Consequently, as $p/q=p-1$, we can obtain that
\[\mathbb{P}(D_n)\leq \inf_{p>1}\big(\mathbb{E}\big[\exp\big\{-(p-1)\big(\lambda x-\lambda^\beta\big)G_n(\beta)\big\}1_{D_n}\big]\big)^{1/p}.\]
The right-hand side of the last inequality attains its maximum at
\[\lambda^{\ast}=(\frac{x}{\beta})^{\frac{1}{\beta-1}}.\]
So
\[\mathbb{P}(D_n)\leq \inf_{p>1}\mathbb{E}\big[\exp\big\{(p-1)(\frac{x}{\beta})^{\frac{\beta}{\beta-1}}(1-\beta)G_n(\beta)\big\}1_{\{S_n\geq xG_n(\beta)\}}\big]^{1/p},\]
which gives the desired inequalities.
\subsection{Proof of Theorem \ref{th2.4}}

The proof of Theorem \ref{th2.4} is similar with the proof of Theorem \ref{th2.2}. Given $1<a<\beta$, introduce the geometric series $b_k=ba^k$ and define random events
\[M_k=\big\{\frac{S_n}{\sqrt[\beta]{G_n(\beta)}}\geq x,(b_k)^{\frac{1}{\beta-1}}\leq \sqrt[\beta]{G_n(\beta)}\leq (b_{k+1})^{\frac{1}{\beta-1}}\big\},\quad k=0,1,\cdots,K,\]
where $K$ stands for the integer part of $\log_aM$. Clearly, it holds
\[\big\{\frac{S_n}{\sqrt[\beta]{G_n(\beta)}}\geq x,(b)^{\frac{1}{\beta-1}}\leq \sqrt[\beta]{G_n(\beta)}\leq (bM)^{\frac{1}{\beta-1}}\big\}\subseteq\cup_{k=0}^{K}M_k,\]
which leads to
\[\mathbb{P}\big(\frac{S_n}{\sqrt[\beta]{G_n(\beta)}}\geq x,(b)^{\frac{1}{\beta-1}}\leq \sqrt[\beta]{G_n(\beta)}\leq (bM)^{\frac{1}{\beta-1}}\big)\leq \sum_{k=0}^{K}\mathbb{P}(M_k).\]
Lemma \ref{t1.4} implies that
\[\mathbb{E}\big[\exp\big\{\lambda S_n-\lambda^\beta G_n(\beta)\big\}1_{M_k}\big]\leq 1.\]
Now, taking $\lambda_k=(\frac{x}{\beta b_k})^{\frac{1}{\beta-1}}$, for any $x>0$, we obtain
\begin{equation*}
\begin{aligned}
1&\geq \mathbb{E}\big[\exp\big\{(\frac{x}{\beta b_k})^{\frac{1}{\beta-1}}S_n-(\frac{x}{\beta b_k})^{\frac{\beta}{\beta-1}}G_n(\beta)\big\}1_{M_k}\big]\\
&\geq \mathbb{E}\big[\exp\big\{(\frac{x}{\beta b_k})^{\frac{1}{\beta-1}}x\sqrt[\beta]{G_n(\beta)}-(\frac{x}{\beta b_k})^{\frac{\beta}{\beta-1}}G_n(\beta)\big\}1_{M_k}\big]\\
&\geq \mathbb{E}\big[\exp\big\{(\frac{1}{\beta b_k})^{\frac{1}{\beta-1}}x^{\frac{\beta}{\beta-1}}\sqrt[\beta]{G_n(\beta)}-(\frac{1}{\beta b_k})^{\frac{\beta}{\beta-1}}x^{\frac{\beta}{\beta-1}}G_n(\beta)\big\}1_{M_k}\big]\\
&\geq \mathbb{E}\big[\exp\big\{\inf_{(b_k)^{\frac{1}{\beta-1}}\leq c\leq (b_{k+1})^{\frac{1}{\beta-1}}}(\frac{1}{\beta b_k})^{\frac{1}{\beta-1}}x^{\frac{\beta}{\beta-1}}c-(\frac{1}{\beta b_k})^{\frac{\beta}{\beta-1}}x^{\frac{\beta}{\beta-1}}c^\beta\big\}1_{M_k}\big]\\
&=\mathbb{E}\big[\exp \big\{x^{\frac{\beta}{\beta-1}}\inf_{(b_k)^{\frac{1}{\beta-1}}\leq c\leq (b_{k+1})^{\frac{1}{\beta-1}}}\big((\frac{1}{\beta b_k})^{\frac{1}{\beta-1}}c-(\frac{1}{\beta b_k})^{\frac{\beta}{\beta-1}}c^\beta\big)\big\}1_{M_k}\big]\\
&=\mathbb{E}\big[\exp \big\{x^{\frac{\beta}{\beta-1}}\inf_{b_k\leq c\leq b_{k+1}}\big((\frac{c}{\beta b_k})^{\frac{1}{\beta-1}}-(\frac{c}{\beta b_k})^{\frac{\beta}{\beta-1}}\big)\big\}1_{M_k}\big]\\
&\geq \mathbb{E}\big[\exp \big\{x^{\frac{\beta}{\beta-1}}\big((\frac{b_{k+1}}{\beta b_k})^{\frac{1}{\beta-1}}-(\frac{b_{k+1}}{\beta b_k})^{\frac{\beta}{\beta-1}}\big)\big\}1_{M_k}\big],
\end{aligned}
\end{equation*}
which implies that
\begin{equation*}
\begin{aligned}
\mathbb{P}(M_k)&\leq \mathbb{E}\big[\exp \big\{-x^{\frac{\beta}{\beta-1}}\big((\frac{a}{\beta})^{\frac{1}{\beta-1}}-(\frac{a}{\beta})^{\frac{\beta}{\beta-1}}\big)\big\}\big]\\
&= \mathbb{E}\big[\exp \big\{-(\frac{x}{\beta})^{\frac{1}{\beta-1}}x\big(a^{\frac{1}{\beta-1}}-\frac{a^{\frac{\beta}{\beta-1}}}{\beta}\big)\big\}\big].
\end{aligned}
\end{equation*}
Let $a=1+(\beta-1)/(1+x)$, then
\[x\big(a^{\frac{1}{\beta-1}}-\frac{a^{\frac{\beta}{\beta-1}}}{\beta}\big)\geq \frac{x}{\beta}\frac{1}{\beta}(\beta-1).\]
So
\begin{equation*}
\begin{aligned}
\mathbb{P}\big(\frac{S_n}{\sqrt[\beta]{G_n(\beta)}}\geq x,&b^{\frac{1}{\beta-1}}\leq \sqrt[\beta]{G_n(\beta)}\leq (bM)^{\frac{1}{\beta-1}}\big)\\
&\leq (1+2(1+x)\ln M)\exp\big\{-(\frac{x}{\beta})^{\frac{\beta}{\beta-1}}(1-\frac{1}{\beta})\big\},
\end{aligned}
\end{equation*}
which gives the desired inequality.
\subsection{Proof of Theorem \ref{th3.2}}
By the least squares estimator of $\theta$,
\[\hat{\theta}_n-\theta=\frac{\sum_{k=1}^{n}\phi_{k-1}X_k}{\sum_{k=1}^{n}\phi_{k-1}^2}.\]
Let
\[\xi_i=\phi_{i-1}\varepsilon_i\qquad \textrm{and}\quad\mathfrak{F}_i=\sigma(\phi_0,\varepsilon_1,\cdots,\phi_{i-1},\varepsilon_i)\quad \textrm{for}\quad i=1,2,\cdots,n.\]
Since we suppose that $\varepsilon_i\leq y$ and $|\phi_i|\leq 1$. The random variable $\varepsilon_{i+1}$ is independent of $\mathfrak{F}_i$, then $(\xi_i,\mathfrak{F}_i)_{i=1,2,\cdots,n}$ is a sequence of martingale differences which satisfies
\[\xi_i\leq y.\]
So
\[\hat{\theta}_n-\theta=\sigma^2\frac{S_n}{\langle S\rangle_n},\]
where
\[\langle S\rangle_n=\sum_{k=1}^{n}\mathbb{E}[\xi_i^2|\mathfrak{F}_{i-1}]=\sigma^2\sum_{k=1}^{n}\phi_{k-1}^2.\]
When $\xi_i\leq y$, then $[S]_n(y)$ vanishes. So $B_n(y)=\langle S\rangle_n$. Applying Theorem \ref{th2.1}, we deduce that for all $x,y>0$,
\begin{equation}\label{eq44}
\begin{aligned}
\mathbb{P}(\hat{\theta}_n-\theta\geq x)&=\mathbb{P}(\frac{S_n}{\langle S\rangle_n}\geq \frac{x}{\sigma^2})\\
&\leq \inf_{p>1}\big(\mathbb{E}\big[\exp\big\{-(p-1)\frac{x^2}{2(\sigma^2+\frac{xy}{3}))}\sum_{k=1}^{n}\phi_{k-1}^2\big\}\big]\big)^{1/p}.
\end{aligned}
\end{equation}
Similarly, we can get that
\begin{equation}\label{eq45}
\mathbb{P}=(\hat{\theta}_n-\theta\leq -x)\leq \inf_{p>1}\big(\mathbb{E}\big[\exp\big\{-(p-1)\frac{x^2}{2(\sigma^2+\frac{xy}{3}))}\sum_{k=1}^{n}\phi_{k-1}^2\big\}\big]\big)^{1/p}.
\end{equation}
Combine (\ref{eq44}) and (\ref{eq45}), we obtain
\begin{equation*}
\mathbb{P}=(|\hat{\theta}_n-\theta|\geq x)\leq 2\inf_{p>1}\big(\mathbb{E}\big[\exp\big\{-(p-1)\frac{x^2}{2(\sigma^2+\frac{xy}{3}))}\sum_{k=1}^{n}\phi_{k-1}^2\big\}\big]\big)^{1/p}.
\end{equation*}
\subsection{Proof of Theorem \ref{th3.3}}
Recall the definition of $\hat{\theta}_n$. It is easy to see that
\[(\hat{\theta}_n-\theta)\sqrt{\sum_{k=1}^{n}\phi_{k-1}^2}=\sigma\frac{S_n}{\sqrt{\langle S\rangle_n}}.\]
Therefore, by Theorem \ref{th2.2}, for all $b>0, M\geq 1$ and $x>0$,
\begin{equation*}
\begin{aligned}
\mathbb{P}&\big((\hat{\theta}_n-\theta)\sqrt{\sum_{k=1}^{n}\phi_{k-1}^2}\geq x,b\leq \sqrt{\sum_{k=1}^{n}\phi_{k-1}^2}\leq bM\big)\\
&\leq \mathbb{P}\big(\frac{S_n}{\sqrt{\langle S\rangle_n}}\geq \frac{x}{\sigma},b\sigma\leq \sqrt{\langle S\rangle_n}\leq \sigma bM\big)\\
&\leq \sqrt{e}\big(1+2(1+\frac{x}{\sigma})\ln M\big)\exp\big\{-\frac{x^2}{2(\sigma^2+xy/3b)}\big\}.
\end{aligned}
\end{equation*}
Similarly,
\begin{equation*}
\begin{aligned}
\mathbb{P}&\big((\hat{\theta}_n-\theta)\sqrt{\sum_{k=1}^{n}\phi_{k-1}^2}\leq -x,b\leq \sqrt{\sum_{k=1}^{n}\phi_{k-1}^2}\leq bM\big)\\
&\leq \sqrt{e}\big(1+2(1+\frac{x}{\sigma})\ln M\big)\exp\big\{-\frac{x^2}{2(\sigma^2+xy/3b)}\big\}.
\end{aligned}
\end{equation*}
Hence, we have for all $b>0, M\geq 1$ and $x>0$,
\begin{equation*}
\begin{aligned}
\mathbb{P}&\big(|\hat{\theta}_n-\theta|\sqrt{\sum_{k=1}^{n}\phi_{k-1}^2}\geq x,b\leq \sqrt{\sum_{k=1}^{n}\phi_{k-1}^2}\leq bM\big)\\
&\leq 2\sqrt{e}\big(1+2(1+\frac{x}{\sigma})\ln M\big)\exp\big\{-\frac{x^2}{2(\sigma^2+ xy/3b)}\big\}.
\end{aligned}
\end{equation*}
\subsection{Proof of Theorem \ref{th3.4}}
Since $T_n$ is $\mathfrak{F}_n$ measurable, by the definition od $d_i$, we can write
\[T_n-\mathbb{E}[T_n]=\sum_{1\leq k\leq n}d_k,\]
where $\{d_k,1\leq k\leq n\}$ is a martingale difference sequence and $d_i\geq 0$ for all $i\geq 1$. So $\langle S\rangle_n(0)$ vanished, then
\[B_n(0)=[S]_n(0)=\sum_{i=1}^{n}d_i^2.\]
By (2.8) in Steele \cite{ste97} and Corollary 5 in Rhee and Talagrand \cite{Rhee87}, we can get that $||d_i||_{\infty}\leq c_1(d)n^{-\frac{1}{d}}$ and $\mathbb{E}[d_i^2]<\infty$ for all $1\leq i\leq n$. By the Cauchy inequality, we can get that
\[\sqrt{\sum_{i=1}^{n}d_i^2}\geq \frac{\sum_{i=1}^{n}d_i}{\sqrt{n}}\geq c_1(d)n^{1/2-1/d}\]
and
\[\sqrt{\sum_{i=1}^{n}d_i^2}\leq c_1(d)\sqrt{\sum_{i=1}^{n}(n-i-1)^{-2/d}}\leq c_1(d)\sqrt{n}.\]
So the condition in Corollary \ref{col2.2} are satisfied, where $b=c_1(d)n^{1/2-1/d}>0$ and $M=n^{1/d}\geq 1$ for all $n\geq 1$.
Then, by inequality (\ref{eq2.7}), we have for all $t>0$,
\begin{equation*}
\begin{aligned}
\mathbb{P}\big(\frac{T_n-\mathbb{E}(T_n)}{\sqrt{\sum_{i=1}^{n}d_i^2}}\geq t, &c_1(d)n^{1/2-1/d}\leq \sqrt{\sum_{i=1}^{n}d_i^2}\leq c_1(d)n^{1/2}\big)\leq\\
&\sqrt{e}(1+\frac{2}{d}(1+t)\ln n)\exp\{-\frac{t^2}{2}\}.
\end{aligned}
\end{equation*}


\end{document}